\documentclass[12pt, a4paper]{article}
\usepackage{amssymb}
\usepackage{latexsym}
\usepackage{amsbsy}

\newcommand{\mput}{\multiput}
\newcommand{\bcen}{\begin{center}}      \newcommand{\ecen}{\end{center}}
\newcommand{\bay}{\begin{array}}      \newcommand{\eay}{\end{array}}
\def\az{\alpha}
\def\bz{\beta}
\def\gz{\gamma}

\def\rank{\mbox{rank }}
\def\min{\mbox{min}}
\def\max{\mbox{max}}
\def\Hom{\mbox{Hom}}

\begin{document}

\title{\Huge {\bf Subrepresentations of Kronecker Representations}}

\author{\Large Yang Han \footnote{The author is supported by Project
10201004 NSFC and OSRF EDC. {\it 2000 Mathematics Subject
Classification}: 16G20, 15A22, 15A21}}

\date{ Institute of Systems Science,
Academy of Mathematics and Systems Science,\\ Chinese Academy of
Sciences, Beijing 100080, P.R. China.\\ E-mail: hany@iss.ac.cn}

\maketitle

\noindent
 {\bf Abstract : } {\footnotesize Translated into the language of
representations of quivers, a challenge in matrix pencil theory is
to find sufficient and necessary conditions for a Kronecker
representation to be a subfactor of another Kronecker representation
in terms of their Kronecker invariants. The problem is reduced to  a
numerical criterion for a Kronecker representation to be a
subrepresentation of another Kronecker representation in terms of
their Kronecker invariants. The key to the problem is the
calculation of ranks of matrices over polynomial rings. For this, a
generalization and specialization approach is introduced. This
approach is applied to provide a numerical criterion for a
preprojective (resp. regular, preinjective) Kronecker representation
to be a subrepresentation of another preprojective (resp. regular,
preinjective) Kronecker representation in terms of their Kronecker
invariants.}

\bigskip

\section*{Introduction}

\medskip

The classification of Kronecker representations was started by
Weierstrass in 1867 and completed by Kronecker in 1890. A natural
problem is to classify the subrepresentations of Kronecker
representations, i.e., the pairs $(N,M)$ in which $N$ is a
subrepresentation of a Kronecker representation $M$, just as has
been done for uniserial rings by Ringel and Schmidmeier \cite{RS}.
However this problem is hopeless to solve completely: Indeed, this
problem is clearly equivalent to classifying those representations
of the quiver
$$\unitlength=1mm
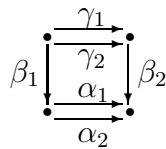
\begin{picture}(20,20) \mput(6,7)(0,10){2}{\circle*{1}}
\mput(17,7)(0,10){2}{\circle*{1}}
\mput(7,6)(0,10){2}{\vector(1,0){9}}
\mput(7,8)(0,10){2}{\vector(1,0){9}}
\mput(6,16)(10.7,0){2}{\vector(0,-1){8}} \put(1,11){$\bz_1$}
\put(18,11){$\bz_2$} \put(10,3){$\az_2$} \put(10,9){$\az_1$}
\put(10,13.5){$\gz_2$} \put(10,19){$\gz_1$}
\end{picture}$$ which satisfy the relations
$\bz_2\gz_1-\az_1\bz_1=\bz_2\gz_2-\az_2\bz_1=0$ and which are such
that the maps $\bz_1$ and $\bz_2$ are inclusion maps \cite{A}. For
this, one would have to classify the representations of the quiver
$$\unitlength=1mm
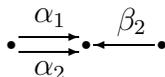
\begin{picture}(20,12) \mput(0,4)(10,0){3}{\circle*{1}}
\mput(1,3)(0,2){2}{\vector(1,0){8}} \put(19,4){\vector(-1,0){8}}
\put(14,6){$\bz_2$} \put(3,7){$\az_1$} \put(3,0){$\az_2$}
\end{picture}$$ which are such that the map $\bz_2$ is an inclusion
map. This problem is clearly wild \cite{R}. Nevertheless we may
study the subrepresentations of Kronecker representations in another
interesting way, namely, to find a numerical criterion for a
Kronecker representation to be a subrepresentation of another
Kronecker representation in terms of their Kronecker invariants.
Later on we will see that the solution of this problem is also the
solution of the first part of the challenge below.

Our original motivation is based on a challenge in matrix pencil
theory. In \cite[p. 329]{LMZZ} the following question, which is
closely related to  pole placement, non-regular feedback, dynamic
feedback, zero placement and early-stage design in control theory is
declared to be a ``challenge'' by the authors.

Recall that a {\it matrix pencil} over a field $\mathbb{K}$ is a
matrix $\lambda E+H$ where $\lambda$ is an indeterminate and $E,H$
are matrices over $\mathbb{K}$ of the same size. Two matrix pencils
$\lambda E_1+H_1$ and $\lambda E_2+H_2$ of the same size are said to
be {\it strictly equivalent}, denoted $\lambda E_1+H_1 \sim \lambda
E_2+H_2$, if there exist invertible constant matrices $P$ and $Q$
such that $\lambda E_1+H_1=P(\lambda E_2+H_2)Q$.

\medskip

{\bf Challenge.} \cite{LMZZ} {\it Let $E, H \in \mathbb{R}^{(m+n)
\times (p+q)}$ and $E', H' \in \mathbb{R}^{m \times p}$. Find
necessary and sufficient conditions in terms of Kronecker invariants
of the matrix pencils $\lambda E+H$ and $\lambda E'+H'$ for the
existence of matrix pencils $F_{12}(\lambda)$, $F_{21}(\lambda)$ and
$F_{22}(\lambda)$ such that $\lambda E+H \sim$
$\left[ \begin{array}{cc} \lambda E'+H'&F_{12}(\lambda)\\
F_{21}(\lambda)&F_{22}(\lambda)
\end{array} \right]$ holds. Moreover, provide an algorithm for
constructing $F_{12}(\lambda)$, $F_{21}(\lambda)$ and
$F_{22}(\lambda)$ whenever a solution exists.}

\medskip

The following  was mentioned in \cite[p. 62]{FS} : ``The problem of
giving necessary and sufficient conditions for the existence of a
matrix pencil with prescribed Kronecker invariants and a prescribed
arbitrary subpencil remains open and seems to be very difficult."
However, partial answers are known when $\lambda E+H$ and $\lambda
E'+H'$ are both regular \cite{B,S,T}; when $\lambda E+H$ is regular
and $\lambda E'+H'$ is arbitrary \cite{CS2}; when $\lambda E+H$ is
arbitrary and $\lambda E'+H'$ is regular \cite{FS}; when $\lambda
E+H$ has rank equal to the number of its rows and $\lambda E'+H'$
has rank equal to the number of its columns \cite{CS1}.

Three approaches, i.e., matrix pencil approach, polynomial approach,
and geometric approach, have been used to attack the Challenge, see
\cite{LMZZ}  and the references cited there. In this paper we
provide the fourth approach, namely representations of quivers. Here
we focus on the first part of the challenge.

The contents of this paper is organized as follows: In section 1, we
first translate the Challenge into the language of representations
of quivers. Thus the Challenge is found equivalent to finding
sufficient and necessary conditions for a Kronecker representation
to be a subfactor of another Kronecker representation in terms of
their Kronecker invariants. Then the problem is reduced to finding a
numerical criterion for a Kronecker representation to be a
subrepresentation of another Kronecker representation in terms of
their Kronecker invariants. And thus the problem becomes fairly
elementary. The key point is to calculate the ranks of matrices over
polynomial rings. Finally we extend the underlying field from the
field of real numbers $\mathbb{R}$ to the field of complex numbers
$\mathbb{C}$ and more generally we work on an arbitrary
algebraically closed field $\mathbb{K}$. Thus the Kronecker
invariants of a Kronecker representation can be expressed simply by
a set of positive integers. In section 2, we consider the
homomorphisms between two Kronecker representations, i.e., the
matrix pairs that satisfy two equations \cite{A}. We partition such
a matrix pair into a block matrix pair. Via easy calculations one
can learn the explicit form of each block in the matrix pair.  This
is very useful. In section 3, we obtain a numerical criterion for a
preprojective (resp. regular, preinjective) Kronecker representation
to be a subrepresentation of another preprojective (resp. regular,
preinjective) Kronecker representation in terms of their Kronecker
invariants. This follows from the calculation of the rank of
matrices over polynomial rings using the generalization and
specialization approach.

\medskip

\section{ Reduction of the Challenge}

\medskip

\subsection{ Translation into the Language of Representations of
Quivers}

Recall that the {\it Kronecker quiver} is the quiver with two
vertices $1, 2$ and two arrows $\az$ and $\bz$ from $1$ to $2$.  A
{\it Kronecker representation} $M$, i.e., a representation of the
Kronecker quiver, can be written as $(M(1),M(2);M(\az),M(\bz))$ or
$(M(\az),M(\bz))$, where $M(1), M(2)$ are the vector spaces
associated with the vertices $1, 2$ respectively and $M(\az),
M(\bz): M(1) \rightarrow M(2)$ are the linear maps that are
represented by the arrows $\az$ and $\bz$, respectively. For more on
representation theory of quivers we refer to \cite{A}. Denote by
${\cal K}$ the representation category of the Kronecker quiver. Note
that in this paper we always consider subrepresentations up to
isomorphism. As a result, we say a Kronecker representation
$N=(N(\az),N(\bz))$ is a {\it subrepresentation} of a Kronecker
representation $M=(M(\az),M(\bz))$ if there is a monomorphism from
$N$ to $M$, or equivalently if there are injective linear maps
$\phi$ and $\psi$ such that $M(\az)\phi=\psi N(\az)$ and
$M(\bz)\phi=\psi N(\bz)$. Dually, a Kronecker representation
$N=(N(\az),N(\bz))$ is called a {\it factor representation} of a
Kronecker representation $M=(M(\az),M(\bz))$ if there is an
epimorphism from $M$ to $N$, or equivalently if there are surjective
linear maps $\phi$ and $\psi$ such that $N(\az)\phi=\psi M(\az)$ and
$N(\bz)\phi=\psi M(\bz)$. A {\it subfactor} of a Kronecker
representation $M$ is a factor representation of a subrepresentation
of $M$, equivalently a subrepresentation of a factor representation
of $M$.

Clearly a matrix pencil $\lambda E+H$ corresponds to a Kronecker
representation $(E,H)$. Moreover, two matrix pencils $\lambda
E_1+H_1$ and $\lambda E_2+H_2$ are strictly equivalent if and only
if $(E_1,H_1)$ and $(E_2,H_2)$ are isomorphic as Kronecker
representations, i.e., if there are invertible matrices $G_1$ and
$G_2$ such that $G_2E_1=E_2G_1$ and $G_2H_1=H_2G_1$. In this way,
the Challenge amounts to {\it finding matrices $E_{12}, E_{21},
E_{22}$, $H_{12}, H_{21}$, and $H_{22}$ such that the two Kronecker
representations $(E,H)$ and $\left(\left[
\begin{array}{cc} E'&E_{12}\\ E_{21}&E_{22} \end{array} \right],
 \left[
\begin{array}{cc} H'&H_{12}\\ H_{21}&H_{22} \end{array}
\right]\right)$ are isomorphic}. If such a solution exists, then we
write $\lambda E+H \succ \lambda E'+H'$ or $\lambda E'+H' \prec
\lambda E+H$ or $(E,H) \succ (E',H')$ or $(E',H') \prec (E,H)$.
Clearly, if $(E_1,H_1) \sim (E,H)$ and $(E,H) \succ (E',H')$, then
$(E_1,H_1) \succ (E',H')$; and if $(E_2,H_2) \sim (E',H')$ and
$(E,H) \succ (E',H')$, then $(E,H) \succ (E_2,H_2)$.

\medskip

{\bf Proposition 1.} {\it $(E,H)\succ(E',H')$ if and only if
$(E',H')$ is a subfactor of $(E,H)$. In particular, the relation
$\succ$ is a partial order on the set of all Kronecker
representations.}

\medskip

{\bf Proof.} If $(E,H)\succ(E',H')$ then there are matrices $E_{12},
E_{21}, E_{22}$, $H_{12}, H_{21}$ and $H_{22}$ such that two
Kronecker representations $(E,H)$ and $\left( \left[
\begin{array}{cc} E'&E_{12}\\ E_{21}&E_{22} \end{array}
\right] , \left[
\begin{array}{cc} H'&H_{12}\\ H_{21}&H_{22} \end{array}
\right]\right)$
 are isomorphic. Since $\left( \left[
\begin{array}{c} I\\ 0 \end{array} \right], I\right)$ is a
monomorphism, $\left( \left[ \begin{array}{c} E'\\
E_{21} \end{array} \right] , \left[ \begin{array}{c} H'\\
H_{21} \end{array} \right]\right)$ is a subrepresentation of
$(E,H)$. Furthermore, $(E',H')$ is a factor representation
of\linebreak
 $\left(\left[ \begin{array}{c} E'\\
E_{21} \end{array} \right] , \left[ \begin{array}{c} H'\\
H_{21} \end{array} \right]\right)$, since $\left(I,  \left[
\begin{array}{cc} I&0 \end{array} \right]\right)$ is an epimorphism.
Thus $(E',H')$ is a subfactor of $(E,H)$. Conversely, if $(E',H')$
is a subfactor of $(E,H)$, then there is a subrepresentation
$(E_1,H_1)$ of $(E,H)$ such that $(E',H')$ is a factor
representation of $(E_1,H_1)$. Hence, there are full rank matrices
$A_i,B_i, i=1,2$, and $A_1', B_2'$ such that $(E,H)A_2=B_2(E_1,H_1),
(E',H')A_1$ $=B_1(E_1,H_1), A_1A_1'$ $=I$, and $B_2'B_2=I$. Since
$B_1B_2'$ and $A_2A_1'$ are full rank matrices, there exist
invertible matrices $C_i,D_i, i=1,2$, such that $B_1B_2'=C_1  \left[
\begin{array}{cc} I&0 \end{array} \right] C_2$ and $A_2A_1'=D_1  \left[
\begin{array}{c} I\\ 0 \end{array} \right] D_2$. Consequently,\\[-7mm]
\begin{eqnarray*}
(E',H')&=& (E',H')A_1A_1'=B_1(E_1,H_1)A_1'
=B_1B_2'B_2(E_1,H_1)A_1'= B_1B_2'(E,H)A_2A_1'\\
&=& C_1 \left[
\begin{array}{cc} I&0
\end{array} \right] C_2(E,H)D_1 \left[ \begin{array}{c} I\\ 0
\end{array} \right] D_2 \sim \left[ \begin{array}{cc} I&0
\end{array} \right] C_2(E,H)D_1  \left[ \begin{array}{c} I\\ 0
\end{array} \right] \\
&\prec& C_2(E,H)D_1 \sim (E,H) \ \ . \hspace*{80mm} \Box
\end{eqnarray*}

\medskip

\subsection{Reduction to the Subrepresentation Case}

Once we find a sufficient and necessary condition $C({\cal N},{\cal
M})$ for a Kronecker representation $N$ to be a subrepresentation of
another Kronecker representation $M$ in terms of the Kronecker
invariants ${\cal N}$ and ${\cal M}$ of $N$ and $M$, then dually we
will find a sufficient and necessary condition $C^*({\cal M},{\cal
N})$ for $N$ to be a factor representation of $M$. Furthermore, we
will find a sufficient and necessary condition for $N$ to be a
subfactor of $M$: There exists a Kronecker module $L$ of Kronecker
invariants ${\cal L}$ such that conditions $C({\cal L},{\cal M})$
and $C^*({\cal L},{\cal N})$ are satisfied. Therefore the question
is reduced from the subfactor one to one of the subrepresentation.

\medskip

{\bf Remark.} The existence question in the condition is not very
easy to handle, but it seems difficult to avoid. Indeed,
 existence question also appear in the results of \cite{CS1,CS2,FS}.

\medskip

\subsection{Extension of the Underlying Field}

Though the question is posed on the field of real numbers
$\mathbb{R}$, we may consider the question on the field of complex
numbers $\mathbb{C}$:

\medskip

{\bf Proposition 2.} {\it A real Kronecker representation $(E',H')$
is a subrepresentation of another real Kronecker representation
$(E,H)$ over $\mathbb{R}$ if and only if the same is the case over
$\mathbb{C}$.}

\medskip

{\bf Proof.} The necessity is trivial. It remains to consider
sufficiency. First, there are full column rank complex matrices $P$
and $Q$ such that $QE'=EP$ and $QH'=HP$. Second, let $P=P_1+iP_2$
and $Q=Q_1+iQ_2$ with $P_j, Q_j, j=1,2,$ being real matrices and
$i=\sqrt{-1}$. Then we have $Q_jE'=EP_j$ and $Q_jH'=HP_j$ for
$j=1,2$. Since $P$ (resp. $Q$) is of full column rank, $P_1=P_2=0$
(resp. $Q_1=Q_2=0$) can not occur. Hence $P_1+xP_2$ (resp.
$Q_1+xQ_2$) is of smaller rank than $P$ (resp. $Q$) for only
finitely many values $x$ in $\mathbb{C}$, i.e., the common roots of
all $\rank P$ (resp. $\rank Q$)--minors of $P_1+xP_2$ (resp.
$Q_1+xQ_2$). Consequently there is some value $x_0$ in $\mathbb{R}$
such that $P_1+x_0P_2$ and $Q_1+x_0Q_2$ are of full column rank, and
$(P_1+x_0P_2,Q_1+x_0Q_2)$ is a monomorphism from $(E',H')$ to
$(E,H)$. \hfill{$\Box$}

\medskip

And more generally, we are able to consider the problem over an
arbitrary algebraically closed field $\mathbb{K}$. By extension of
the underlying field we can simply express Kronecker invariants as a
set of integers (see section 2.1 below), this is of great benefit.

\medskip

\section{ Homomorphisms between two Kronecker
Representations}

\medskip

\noindent Note that a homomorphism between two Kronecker
representations is just a pair of matrices satisfying two equations.
In this section we partition these two matrices in the natural way
(corresponding to their direct sum decompositions of indecomposable
representations) and observe the form of every block.

\medskip

\subsection{Kronecker Invariants}

Denote by $I$ the identity matrix and by $J$ the Jordan block with
eigenvalue 0 (of the appropriate size). Denote by
$\mathbb{P}^1(\mathbb{K})$ the projective line over $\mathbb{K}$. By
the well-known Krull-Schmidt theorem, a Kronecker representation can
be decomposed into a direct sum of indecomposable Kronecker
representations. Let $Q_i:=(\mathbb{K}^{i-1},\mathbb{K}^i; {\tiny
\left[\begin{array}{c} I\\0
\end{array} \right]}, {\tiny \left[\begin{array}{c} 0\\I
\end{array} \right]})$,
$R_{\infty,i}:= (\mathbb{K}^i,\mathbb{K}^i;J,I), R_{p,i}:=
(\mathbb{K}^i,\mathbb{K}^i;I,pI+J)$, and
$J_i:=(\mathbb{K}^i,\mathbb{K}^{i-1}; {\tiny \left[\begin{array}{cc}
I&0
\end{array} \right]},$ ${\tiny \left[\begin{array}{cc} 0&I
\end{array} \right]}), p \in \mathbb{K}, i \in \mathbb{N}_1:=\{1,2,...\}$. Then
the sets $\{Q_i|i \in \mathbb{N}_1\}, \{R_{p,i}|p \in
\mathbb{P}^1(\mathbb{K}), i \in \mathbb{N}_1\}$ and $\{J_i|i \in
\mathbb{N}_1\}$, called {\it preprojective, regular,} and {\it
preinjective} indecomposable Kronecker representations respectively,
constitute a complete set of nonisomorphic indecomposable Kronecker
representations \cite{A}. Up to isomorphism, a Kronecker
representation $M$ can be uniquely written as $M=(\oplus
^{m^P}_{i=1}Q_{a_i}) \oplus (\oplus_{p \in \mathbb{P}^1(\mathbb{K})}
\oplus_{i=1}^{m^p} R_{p,b^p_i}) \oplus (\oplus ^{m^I}_{i=1}J_{c_i})$
for some positive integers $a_i, i=1,...,m^P; b^p_i, i=1,...,m^p, p
\in \mathbb{P}^1(\mathbb{K}); c_i, i=1,...,m^I$ (notice that the
superscripts do not mean power). The Kronecker representation $M$ is
uniquely determined by $a_i, b^p_i, c_i$, which are called the {\it
Kronecker invariants} of $M$. Moreover, a Kronecker representation
is said to be {\it preprojective} (resp. {\it regular,
preinjective}) if it is the direct sum of preprojective (resp.
regular, preinjective) indecomposable representations.

\medskip

{\bf Remark.} Usually the Kronecker invariants of $M$ viewed as a
matrix pencil are referred to the {\it row minimal indices}, the
{\it infinite elementary factors}, the {\it finite elementary
factors}, and the {\it column minimal indices} \cite{G,LMZZ}. Over
an algebraically closed field $\mathbb{K}$, they correspond to
positive integers $a_i, b^{\infty}_i, b^p_i (p \in \mathbb{K}), c_i$
respectively.

\medskip

\subsection{Decomposition of Homomorphism}

Let $M$ and $N$ be two Kronecker representations. Then $M=M^P \oplus
M^R \oplus M^I$ and $N=N^P \oplus N^R \oplus N^I$ where $M^P=\oplus
^{m^P}_{i=1}Q_{a_i}$ with $a_1 \geq a_2 \geq \cdots \geq a_{m^P}$,
$M^R=\oplus_{p \in \mathbb{P}^1(\mathbb{K})} \oplus_{i=1}^{m^p}
R_{p,b^p_i}$ with $b^p_1 \geq b^p_2 \geq \cdots \geq b^p_{m^p}$ for
every $p \in \mathbb{P}^1(\mathbb{K})$, $M^I=\oplus
^{m^I}_{i=1}J_{c_i}$ with $c_1 \geq c_2 \geq \cdots \geq c_{m^I}$,
$N^P=\oplus ^{n^P}_{i=1}Q_{d_i}$ with $d_1 \geq d_2 \geq \cdots \geq
d_{n^P}$, $N^R=\oplus_{p \in \mathbb{P}^1(\mathbb{K})}
\oplus_{i=1}^{n^p} R_{p,e^p_i}$ with $e^p_1 \geq e^p_2 \geq \cdots
\geq e^p_{n^p}$ for every $p \in \mathbb{P}^1(\mathbb{K})$,
$N^I=\oplus ^{n^I}_{i=1}J_{f_i}$ with $f_1 \geq f_2 \geq \cdots \geq
f_{n^I}$. Of course these numbers $a_i,b^p_i,c_i,d_i,e^p_i,f_i$ are
positive integers. Once again the superscripts do not mean power
here.

\medskip

By [1; Theorem 7.5], any homomorphism of representations $\phi \in
\Hom_{\cal K}(N,$ $M)$ can be written as $\phi= {\tiny
\left[\begin{array}{ccc} \phi^{PP}&0&0\\ \phi^{RP}&\phi^{RR}&0\\
\phi^{IP}&\phi^{IR}&\phi^{II} \end{array} \right]}$ where $\phi^{ST}
\in \Hom_{\cal K}(N^T,$ $M^S)$ for $S,T \in \{P,R,I\}$. Writing
$\phi$ as a matrix pair, we have
$$\phi=(\phi^1,\phi^2)=
\left({\tiny
\left[\begin{array}{ccc} \phi^{PP1}&0&0\\ \phi^{RP1}&\phi^{RR1}&0\\
\phi^{IP1}&\phi^{IR1}&\phi^{II1} \end{array} \right]},{\tiny
\left[\begin{array}{ccc} \phi^{PP2}&0&0\\ \phi^{RP2}&\phi^{RR2}&0\\
\phi^{IP2}&\phi^{IR2}&\phi^{II2} \end{array} \right]}\right)$$ with
$\phi^{ST}=(\phi^{ST1},\phi^{ST2})$  for $S,T \in \{P,R,I\}$.

\medskip

\subsection{Analysis of $\phi^{PP}$ and $\phi^{II}$}

We can write $M^P=(M^P(1),M^P(2); $  $M^P(\az),$ $M^P(\bz))$ and
$N^P=(N^P(1), N^P(2);$\linebreak $ N^P(\az), N^P(\bz))$, where
$M^P(1)=\mathbb{K}^{\sum _{i=1}^{m^P} (a_i-1)}$,
$M^P(2)=\mathbb{K}^{\sum _{i=1}^{m^P} a_i}$, $N^P(1)=
\mathbb{K}^{\sum_{j=1}^{n^P} (d_j-1)}$, $N^P(2)=\mathbb{K}^{\sum
_{j=1}^{n^P} d_j}$, $M^P(\az)$ and $N^P(\az)$ are of the form diag
$\{{\tiny \left[\begin{array}{c} I\\0 \end{array} \right]}, \ldots ,
{\tiny \left[\begin{array}{c} I\\0 \end{array} \right]}\}$, and
$M^P(\bz)$ and $N^P(\bz)$ are of the form diag$\{{\tiny
\left[\begin{array}{c} 0\\I \end{array} \right]}, \ldots , {\tiny
\left[\begin{array}{c} 0\\I \end{array} \right]}\}$. We can write
$\phi^{PP}=(\phi^{PP1},\phi^{PP2})$ where $\phi^{PP1}$ and
$\phi^{PP2}$ are $(\sum _{i=1}^{m^P} (a_i-1)) \times
(\sum_{j=1}^{n^P} (d_j-1))$ and $(\sum _{i=1}^{m^P} a_i) \times
(\sum_{j=1}^{n^P} d_j)$ matrices respectively. By partitioning into
$m^P \times n^P$ block matrices in the natural way (corresponding to
their direct sum decomposition), we have
$\phi^{PP1}=(\phi^{PP1}_{ij})_{ij}$ and
$\phi^{PP2}=(\phi^{PP2}_{ij})_{ij}, i=1,...,m^P, j=1,...,n^P$. Since
$M^P(\az)\phi^{PP1}=\phi^{PP2}N^P(\az)$ and
$M^P(\bz)\phi^{PP1}=\phi^{PP2}N^P(\bz)$, we have ${\tiny
\left[\begin{array}{c} I\\0 \end{array} \right]} \phi^{PP1}_{ij}=
\phi^{PP2}_{ij} {\tiny \left[\begin{array}{c} I\\0 \end{array}
\right]}$ and ${\tiny \left[\begin{array}{c} 0\\I
\end{array} \right]} \phi^{PP1}_{ij}= \phi^{PP2}_{ij} {\tiny
\left[\begin{array}{c} 0\\I \end{array} \right]}$. Therefore the
blocks $\phi^{PP1}_{ij}$ and $\phi^{PP2}_{ij}$ have the form
\begin{equation}\label{star}
\left[\begin{array}{ccc} x^{PPij}_1&&\\ \ddots&\ddots&\\
x^{PPij}_{a_i-d_j+1}&\ddots&x^{PPij}_1\\
&\ddots&\ddots\\ &&x^{PPij}_{a_i-d_j+1} \end{array} \right]
\end{equation}
of size $(a_i-1) \times (d_j-1)$ and $a_i \times d_j$, respectively,
in case $a_i \geq d_j$, and empty  otherwise.

\medskip

Similarly we can write $\phi^{II}=(\phi^{II1},\phi^{II2})$ where
$\phi^{II1}$ and $\phi^{II2}$ are $(\sum _{i=1}^{m^I} c_i) \times
(\sum_{j=1}^{n^I} f_j)$ and $(\sum _{i=1}^{m^I} (c_i-1)) \times
(\sum_{j=1}^{n^I} (f_j-1))$ matrices, respectively. We partition
these into $m^I \times n^I$ block matrices in the natural way and
have $\phi^{II1}=(\phi^{II1}_{ij})_{ij}$ and
$\phi^{II2}=(\phi^{II2}_{ij})_{ij}$ where the blocks
$\phi^{II1}_{ij}$ and $\phi^{II2}_{ij}$ have the form
\begin{equation}\label{2star}
\left[\begin{array}{ccccc} x^{IIij}_1&\ddots&x^{IIij}_{f_j-c_i+1}&&\\
&\ddots&\ddots&\ddots&\\ &&x^{IIij}_1&\ddots&x^{IIij}_{f_j-c_i+1}
\end{array} \right]
\end{equation}
 of size $c_i \times f_j$ and $(c_i-1) \times
(f_j-1)$, respectively, in case $c_i \leq f_j$, and empty otherwise.

\medskip

\subsection{Analysis of $\phi^{RR}$}

Note that we can write $M^R=(M^R(1),M^R(2); M^R(\az),M^R(\bz))$ and
$N^R=(N^R(1),$\linebreak $ N^R(2); N^R(\az),N^R(\bz))$, where
$M^R(1)=M^R(2)=\mathbb{K}^{\sum _{p \in \mathbb{P}^1(\mathbb{K})}
\sum_{i=1}^{m^p} b^p_i}$, $N^R(1)=N^R(2)=\mathbb{K}^{\sum _{p \in
\mathbb{P}^1(\mathbb{K})} \sum_{j=1}^{n^p} e^p_j}$, $M^R(\az)$ and
$N^R(\az)$ are of the form diag$\{J, \ldots , J, I, \ldots, I \}$,
and $M^R(\bz)$ and $N^R(\bz)$ are of the form diag$\{I, \ldots, I,
\ldots, pI+J, \ldots, pI+J, \ldots \}$. We can write
$\phi^{RR}=(\phi^{RR1},\phi^{RR2})$ where $\phi^{RR1}$ and
$\phi^{RR2}$ are $(\sum _{p \in \mathbb{P}^1(\mathbb{K})}
\sum_{i=1}^{m^p} b^p_i) \times (\sum _{p \in
\mathbb{P}^1(\mathbb{K})} \sum_{j=1}^{n^p} e^p_j)$ matrices. By [1;
Theorem 7.5], we have $\phi^{RR1}=$ diag$\{\phi^{RR1p}\}_{p \in
\mathbb{P}^1(\mathbb{K})}$ and $\phi^{RR2}=$ diag$\{\phi^{RR2p}\}_{p
\in \mathbb{P}^1(\mathbb{K})}$. If we partition these matrices into
$(\sum _{p \in \mathbb{P}^1(\mathbb{K})}m^p) \times (\sum _{p \in
\mathbb{P}^1(\mathbb{K})}n^p)$ block matrices in the natural way, we
have $\phi^{RR1p}=(\phi^{RR1p}_{ij})_{ij}$ and
$\phi^{RR2p}=(\phi^{RR2p}_{ij})_{ij}$. Since
$M^R(\az)\phi^{RR1}=\phi^{RR2}N^R(\az)$ and
$M^R(\bz)\phi^{RR1}=\phi^{RR2}N^R(\bz)$, we have $J \phi^{RR1
\infty}_{ij}=\phi^{RR2 \infty}_{ij} J$ and $I \phi^{RR1
\infty}_{ij}= \phi^{RR2 \infty}_{ij} I$, $I
\phi^{RR1p}_{ij}=\phi^{RR2p}_{ij} I$, and $(pI+J) \phi^{RR1p}_{ij}=
\phi^{RR2p}_{ij} (pI+J)$ for every $p \in \mathbb{K}$. Therefore the
block $\phi^{RR1p}_{ij}=\phi^{RR2p}_{ij}$ has  the form
\begin{equation}\label{3star}
\left[\begin{array}{ccccc} &x^{ppij}_{g^p_{ij}}&\ddots&
x^{ppij}_2&x^{ppij}_1\\ &&\ddots&\ddots&x^{ppij}_2\\
&&&\ddots&\ddots\\ &&&&x^{ppij}_{g^p_{ij}}\\ &&&& \end{array}
\right]
\end{equation}
of size $b^p_i \times e^p_j$, where $g^p_{ij}=\min \{b^p_i, e^p_j\}$
for every $p \in \mathbb{P}^1(\mathbb{K})$ and $h \in \{1,2\}$.

\medskip

{\bf Remark.} In a similar way, one can easily describe every block
in $\phi^1$ and $\phi^2$. However, as sections 2.3 and 2.4 are
enough for later use, all other cases are omitted here.

\medskip

\section{ Subrepresentations of Kronecker Representations}

\medskip

Assume that $N'$ and $M'$ are preprojective (resp. regular, or
preinjective) Kronecker representations. In this section we provide
a sufficient and necessary condition for $N'$ to be a
subrepresentation of $M'$ in terms of their Kronecker invariants.
For convenience, we consider the preprojective (resp. regular, or
preinjective) parts of the Kronecker modules $N$ and $M$ given in
section 2.2 instead of $N'$ and $M'$.

A {\it generic matrix} is a matrix whose elements are pairwise
different indeterminates. A matrix pair $\phi=(\phi^1,\phi^2)$ is
called a {\it generic homomorphism} from $N$ to $M$ if $\phi^1$ and
$\phi^2$ are generic matrices satisfying $M(\az)\phi^1=\phi^2N(\az)$
and $M(\bz)\phi^1=\phi^2N(\bz)$. Clearly, a generic homomorphism
from $N$ to $M$ is a homomorphism from $N$ to $M$ over some
transcendental extension field of $\mathbb{K}$. Once the
indeterminates in the generic homomorphism $\phi$ take special
values in $\mathbb{K}$ then $\phi$ becomes a homomorphism from $N$
to $M$. Conversely, any homomorphism from $N$ to $M$ can be obtained
in this way. From now on $\phi=(\phi^1,\phi^2)$ is always assumed to
be a generic homomorphism from $N$ to $M$. Clearly, $N$ is a
subrepresentation of $M$ if and only if there exists a monomorphism
from $N$ to $M$, or if and only if the generic homomorphism
$\phi=(\phi^1,\phi^2)$ from $N$ to $M$ is a monomorphism over some
rational function field over $\mathbb{K}$, or if and only if
$\phi^1$ and $\phi^2$ viewed as matrices over the polynomial rings,
equivalently over their quotient fields, are of full column ranks.
If we partition $\phi$ as  done in section 2.2 and partition
$\phi^{SSh}, S\in \{P,R,I\}, h \in \{1,2\}$, as  done in section 2.3
and section 2.4, then the blocks in $\phi^{SSh}$ have the forms
(\ref{star}), (\ref{2star}) or (\ref{3star}) in sections 2.3 and 2.4
where all $x^*_*$ are assumed to be indeterminates. Thus $N^S$ is a
subrepresentation of $M^S$ if and only if $\phi^{SS1}$ and
$\phi^{SS2}$ are of full column rank. In order to determine when
$\phi^{SS1}$ and $\phi^{SS2}$ are of full column rank, we calculate
the ranks of $\phi^{SS1}$ and $\phi^{SS2}$.

\medskip

\subsection{Generalization and Specialization}

In order to calculate the ranks of the matrices $\phi^{SS1}$ and
$\phi^{SS2}$, we employ the {\it generalization and specialization
approach}. The {\it generalization} procedure consists of replacing
some elements in the matrix of rational functions $A$ with new
independent indeterminates, so that the rank of the resulting matrix
of rational functions provides an upper bound for the rank of the
original matrix $A$. The {\it specialization} procedure consists of
replacing some indeterminates in $A$ with special values, usually 0
or 1, so that the rank of the resulting matrix provides a lower
bound for the rank of the original matrix $A$. Usually, by a series
elementary transformations of matrices and generalizations, we can
obtain a matrix of rational functions $B$ from $A$, and by
specialization we can obtain a matrix $C$ from $A$. It will be shown
that $\rank B=\rank C$. Thus we conclude that
 $\rank
A=\rank B=\rank C$. In the following we will apply this approach to
calculate the ranks of $\phi^{SS1}$ and $\phi^{SS2}$.

\medskip

First we calculate the ranks of block upper triangular generic
matrices by the generalization-specialization approach. The rank
formula obtained is closely related to the rank formula obtained in
the preprojective-to-preprojective and preinjective-to-preinjective
cases (see the remarks in section 3.3 and section 3.4 below).

\medskip

{\bf Proposition 3.} {\it Let $A=(A_{ij})_{ij}$ with $1 \leq i,j
\leq q$ be a block upper triangular generic matrix, i.e., $A_{ij}=0$
for $1 \leq j <i \leq q$ and $A_{ij}$ is $r_i \times c_j$ generic
matrix for $1 \leq i \leq j \leq q$. Assume that all indeterminates
in $A$ are different. Then} $\rank  A= \min
\{\sum^i_{j=1}r_j+\sum^q_{j=i+1}c_j|0 \leq i \leq q\}$.

\medskip

{\bf Remark.} By convention we require $\sum_{j=k}^iy_j=0$ if $i<k$.

\medskip

{\bf Proof.} Let $F$ be the transcendental extension field of
$\mathbb{K}$ obtained by adding all indeterminates in $A$, i.e., the
field of rational functions in all indeterminates in $A$ over
$\mathbb{K}$. We proceed by induction on $q$: It is trivial for case
$q=1$. Now consider the case $q \geq 2$.

If $r_1 \leq c_1$ then by elementary transformations over $F$, $A$
can be reduced to another block upper triangular matrix
$A'=(A'_{ij})_{ij}$, $1 \leq i,j \leq q$, where $A'_{ij}=0$ for $1
\leq j <i \leq q$, $A'_{11}=[I,0]$ and $I$ is the $r_1 \times r_1$
identity matrix, and $A'_{1j}=0$ for $2 \leq j \leq q$ with
$A'_{ij}=A_{ij}$ for $2 \leq i \leq q$. By our induction hypothesis,
$\rank  A=r_1 + \min \{\sum^i_{j=2}r_j+\sum^q_{j=i+1}c_j|1 \leq i
\leq q\} = \min \{\sum^i_{j=1}r_j+\sum^q_{j=i+1}c_j|0 \leq i \leq
q\}$.

If $r_1>c_1$ then by elementary transformations over $F$, $A$ can be
reduced to another block upper triangular matrix
$A'=(A'_{ij})_{ij}$, $1 \leq i,j \leq q$, where $A'_{ij}=0$ for $1
\leq j <i \leq q$, $A'_{11}=I$ for the $c_1 \times c_1$ identity
matrix, and $A'_{1j}=0$ for $2 \leq j \leq q$ with $A'_{ij}=A_{ij}$
for $3 \leq i \leq q$. By generalization, i.e., replacing all
elements in the $(r_1+r_2-c_1) \times c_j$ matrices $A'_{2j}$, $2
\leq j \leq q,$ with different new indeterminates, we obtain a
matrix $B$. By induction hypothesis, we get $\rank  A= \rank  A'
\leq \rank  B = c_1 + \min \{\sum^q_{j=2}c_j,
(r_1+r_2-c_1)+\sum^q_{j=3}c_j, (r_1+r_2-c_1)+r_3+\sum^q_{j=4}c_j,
... , (r_1+r_2-c_1)+\sum^q_{j=3}r_j \}= \min
\{\sum^i_{j=1}r_j+\sum^q_{j=i+1}c_j|0 \leq i \leq q\}$. On the other
hand, by specialization, i.e., taking the $(1,1), (2,2), \cdots
,(c_1,c_1)$ entries of $A$ to be 1 and all other indeterminates
lying in the same rows or columns as these entries as 0. The
resulting matrix $C$  clearly has the same rank as $B$. Thus $\rank
A \geq \rank C =\rank B$. Finally $\rank A = \rank C =\rank B= \min
\{\sum^i_{j=1}r_j+\sum^q_{j=i+1}c_j|0 \leq i \leq q\}$.
\hfill{$\Box$}

\medskip

\subsection{The Preprojective to Preprojective Case}

Keeping in mind the analysis of $\phi^{PP}$ in section 2.3, let
$$\begin{array}{l} r_1:=\max \{1 \leq j \leq
n^P|d_j>a_1\};\\ s_1:=\max \{1\leq i \leq m^P|d_{r_1+1} \leq a_i\};\\
...\\ r_l:=\max \{1 \leq j \leq n^P|d_j>a_{s_{l-1}+1}\};\\
s_l:=\max \{1\leq i \leq m^P|d_{r_l+1} \leq a_i\};\\ ...\\
r_t=n^P.
\end{array}$$ Note that $r_1$ is just the number of zero blocks in the
first block row of $\phi^{PP2}$, $s_1$ is just the number of the
block rows of $\phi^{PP2}$ having the largest number of nonzero
blocks. In the following all undefined numbers such as $s_0$ are
assumed to be 0.

\medskip

{\bf Proposition 4.} $\rank  \phi^{PP2}= \min \{\sum^{s_i}_{j=1}a_j+
\sum^{r_t}_{j=r_{i+1}+1}d_j|0 \leq i \leq t-1\}.$

\medskip

{\bf Proof.} We calculate $\rank  \phi^{PP2}$ by induction on $t$.
If $t=1$ then $\phi^{PP2}= 0$ and we are done. Assume $t \geq 2$.

{\it Case 1.} $\sum^{s_i}_{j=1}a_j \geq \sum^{r_{i+1}}_{j=r_1+1}d_j,
1 \leq i \leq t-2$.

In this case we do not need to use induction. Clearly
$$\rank
\phi^{PP2} \leq \min \{\sum^{s_{t-1}}_{j=1}a_j,
\sum^{r_t}_{j=r_1+1}d_j\}= \min \{\sum^{s_i}_{j=1}a_j+
\sum^{r_t}_{j=r_{i+1}+1}d_j|0 \leq i \leq t-1\}\ .$$
 Next we prove that
$\rank  \phi^{PP2} \geq \min \{\sum^{s_{t-1}}_{j=1}a_j,
\sum^{r_t}_{j=r_1+1}d_j\}$. We proceed by specialization, namely we
let the indeterminates in $\phi^{PP2}$ take the special values 0 or
1 such that the resulting matrix is of rank $\min
\{\sum^{s_{t-1}}_{j=1}a_j, \sum^{r_t}_{j=r_1+1}d_j\}$.

(1) If $a_1 \geq \sum^{r_t}_{j=r_1+1}d_j$ then let the
$(1,\sum^{r_1}_{j=1}d_j+1), (2,\sum^{r_1}_{j=1}d_j+2), ..., $
\linebreak $ (\sum^{r_t}_{j=r_1+1}d_j, \sum^{r_t}_{j=1}d_j)$
elements of $\phi^{PP2}$ take 1, and let all other indeterminates
take 0. This finishes the specialization.

(2) If $\sum^{u_1}_{j=r_1+1}d_j \leq a_1 <
\sum^{u_1+1}_{j=r_1+1}d_j$ for some $r_1+1 \leq u_1 < r_t$ then we
let the $(1,\sum^{u_1}_{j=1}d_j+1), (2,\sum^{u_1}_{j=1}d_j+2), ...,
(a_1-\sum^{u_1}_{j=r_1+1}d_j, a_1+\sum^{r_1}_{j=1}d_j),
(a_1-\sum^{u_1}_{j=r_1+1}d_j+1, \sum^{r_1}_{j=1}d_j+1),
(a_1-\sum^{u_1}_{j=r_1+1}d_j+2, \sum^{r_1}_{j=1}d_j+2), ...,
(a_1,\sum^{u_1}_{j=1}d_j)$ entries of $\phi^{PP2}$ take the value 1.
If $t=2$ then  all other indeterminates are set to 0. This ends the
specialization.

(3) For $t \geq 3$ suppose $a_2<d_{u_1+1}$. Then $s_1=1$ and $r_2
\geq u_1+1$. This contradicts  the assumption $\sum^{s_1}_{j=1}a_j
\geq \sum^{r_2}_{j=r_1+1}d_j$. Thus $a_2 \geq d_{u_1+1}$. If
$a_1+a_2 \geq \sum^{r_t}_{j=r_1+1}d_j$ then we set the
$(a_1+a_2-(\sum^{u_1+1}_{j=r_1+1}d_j-a_1)+1,
\sum^{r_1}_{j=1}d_j+a_1+1),
(a_1+a_2-(\sum^{u_1+1}_{j=r_1+1}d_j-a_1)+2,
\sum^{r_1}_{j=1}d_j+a_1+2), ..., (a_1+a_2,\sum^{u_1+1}_{j=1}d_j);
(a_1+1,\sum^{u_1+1}_{j=1}d_j+1), (a_1+2,\sum^{u_1+1}_{j=1}d_j+2),
..., (a_1+\sum^{r_t}_{j=u_1+2}d_j,\sum^{r_t}_{j=1}d_j)$ entries of
$\phi^{PP2}$ equal to 1, and choose all other indeterminates as 0.

(4) If $\sum^{u_2}_{j=r_1+1}d_j \leq a_1+a_2<
\sum^{u_2+1}_{j=r_1+1}d_j$ for some $r_1+1 \leq u_2 <r_t$, then we
set the
$(a_1+a_2-(\sum^{u_1+1}_{j=r_1+1}d_j-a_1)+1,\sum^{r_1}_{j=1}d_j+a_1+1),
(a_1+a_2-(\sum^{u_1+1}_{j=r_1+1}d_j-a_1)+2,\sum^{r_1}_{j=1}d_j+a_1+2),...,
(a_1+a_2,\sum^{u_1+1}_{j=1}d_j);(a_1+1,\sum^{u_2}_{j=1}d_j+1),
(a_1+2,\sum^{u_2}_{j=1}d_j+2),...,
(a_1+(a_1+a_2-\sum^{u_2}_{j=r_1+1}d_j),
\sum^{r_1}_{j=1}d_j+a_1+a_2);
(a_1+(a_1+a_2-\sum^{u_2}_{j=r_1+1}d_j)+1,\sum^{u_1+1}_{j=1}d_j+1),
(a_1+(a_1+a_2-\sum^{u_2}_{j=r_1+1}d_j)+2,\sum^{u_1+1}_{j=1}d_j+2),...,
(a_1+a_2-(\sum^{u_1+1}_{j=r_1+1}d_j-a_1),\sum^{u_2}_{j=1}d_j)$
elements of $\phi^{PP2}$ equal to 1. If $t=3$ we set all other
indeterminates equal to 0.

(5) For $t \geq 4$ suppose $a_3<d_{u_2+1}$. Then there exists some
$s_i = 2$ with $1 \leq i \leq 2$ such that  $r_{i+1} \geq u_2+1$.
This contradicts  the assumption $\sum^{s_i}_{j=1}a_j \geq
\sum^{r_{i+1}}_{j=r_1}d_j$. Thus $a_3 \geq d_{u_2+1}$. If
$a_1+a_2+a_3 \geq \sum^{r_t}_{j=r_1+1}d_j$, then we let the
$(a_1+a_2+a_3-(\sum^{u_2+1}_{j=r_1+1}d_j-a_1-a_2)+1,
\sum^{r_1}_{j=1}d_j+a_1+a_2+1),
(a_1+a_2+a_3-(\sum^{u_2+1}_{j=r_1+1}d_j-a_1-a_2)+2,
\sum^{r_1}_{j=1}d_j+a_1+a_2+2), ...,
(a_1+a_2+a_3,\sum^{u_2+1}_{j=1}d_j),(a_1+a_2+1,\sum^{u_2+1}_{j=1}d_j+1);
(a_1+a_2+2,\sum^{u_2+1}_{j=1}d_j+2),...,
(a_1+a_2+\sum^{r_t}_{j=u_2+2}d_j,\sum^{r_t}_{j=1}d_j)$ elements of
$\phi^{PP2}$ be 1, and set all other indeterminates to be 0.

Proceeding in this way, this process will end with one of two
possibilities:

(i) we can proceed in $2s_{t-1}$ steps: In this case all nonzero
rows are exhausted.

(ii) we can proceed in $2q-1$ steps with $1 \leq q \leq s_{t-1}$: In
this case all nonzero columns are exhausted.

Via our specializations we have obtained a (0,1)-matrix whose rank
is $\sum^{s_{t-1}}_{j=1}a_j$ (resp. $\sum^{r_t}_{j=r_1+1}d_j$) in
the case (i) (resp. (ii)): Indeed this (0,1)-matrix can be reduced
by elementary transformations to a (0,1)-matrix for which in  case
(i) (resp. (ii)) there are just $\sum^{s_{t-1}}_{j=1}a_j$ (resp.
$\sum^{r_t}_{j=r_1+1}d_j$) elements  1 lying in different rows and
columns (by keeping the ones as far to the left as possible).

{\it Case 2.} Assume  that
$\sum^{s_i}_{j=1}a_j<\sum^{r_{i+1}}_{j=r_1+1}d_j$ for some $1 \leq i
\leq t-2$ and let $v:=\min \{1 \leq i \leq t-2|
\sum^{s_i}_{j=1}a_j<\sum^{r_{i+1}}_{j=r_1+1}d_j\}$. Let
$A:=\phi^{PP2}(1,...,\sum^{s_v}_{j=1}a_j;
\sum^{r_1}_{j=1}d_j+1,...,\sum^{r_{v+1}}_{j=1}d_j)$ be the submatrix
of $\phi^{PP2}$ which is the intersection of the
$1$-st,...,$(\sum^{s_v}_{j=1}a_j)$-th rows of $\phi^{PP2}$ and the
$(\sum^{r_1}_{j=1}d_j+1)$-st,...,$(\sum^{r_{v+1}}_{j=1}d_j)$-th
columns of $\phi^{PP2}$. By case 1, we have $\rank A=
\sum^{s_v}_{j=1}a_j$. By the induction hypothesis, the rank of the
submatrix
$B:=\phi^{PP2}(\sum^{s_v}_{j=1}a_j+1,...,\sum^{s_{t-1}}_{j=1}a_j;
\sum^{r_{v+1}}_{j=1}d_j+1,...,\sum^{r_t}_{j=1}d_j)$ of $\phi^{PP2}$
is equal to $\min
\{\sum^{s_i}_{j=s_v+1}a_j+\sum^{r_t}_{j=r_{i+1}+1}d_j|v\leq i \leq
t-1\}$. Thus $\rank  \phi^{PP2} = \rank  A + \rank  B =
\sum^{s_v}_{j=1}a_j + \min \{\sum^{s_i}_{j=s_v+1}a_j+
\sum^{r_t}_{j=r_{i+1}+1}d_j|v\leq i \leq t-1\} = \min
\{\sum^{s_i}_{j=1}a_j+ \sum^{r_t}_{j=r_{i+1}+1}d_j|0\leq i \leq
t-1\}$. \hfill{$\Box$}

\medskip

Note that $r_1=\max \{1 \leq j \leq n^P|d_j-1>a_1-1\}$; $s_1=\max
\{1\leq i \leq m^P|d_{r_1+1}-1 \leq a_i-1\}$; ... ; $r_l=\max \{1
\leq j \leq n^P|d_j-1>a_{s_{l-1}+1}-1\}$; $s_l=\max \{1\leq i \leq
m^P|d_{r_l+1}-1 \leq a_i-1\}$; ... ; $r_t=n^P$. By Proposition 4 we
have the following formula on $\rank \phi^{PP1}$.

\medskip

{\bf Corollary 5.} $\rank \phi^{PP1}=\min \{\sum^{s_i}_{j=1}(a_j-1)+
\sum^{r_t}_{j=r_{i+1}+1}(d_j-1)|0 \leq i \leq t-1\}.$

\medskip

{\bf Remark.} By Proposition 4, Corollary 5, and Proposition 3 we
find that $\rank  \phi^{PPh}$, $h \in \{1,2\}$, is equal to the rank
of the matrix obtained from $\phi^{PPh}$ by replacing each nonzero
block in $\phi^{PPh}$ with a generic matrix of the same size. (Of
course all indeterminates in these generic matrices are assumed to
be different.)

\medskip

By Proposition 4 and Corollary 5 we obtain a numerical criterion for
a preprojective Kronecker representation to be a subrepresentation
of another preprojective Kronecker representation in terms of their
Kronecker invariants.

\medskip

{\bf Theorem 6.} {\it $N^P$ is a subrepresentation of $M^P$ if and
only if $r_1=0,$ $\sum^{s_i}_{j=1}a_j \geq
\sum^{r_{i+1}}_{j=r_1+1}d_j$ and $\sum^{s_i}_{j=1}(a_j-1) \geq
\sum^{r_{i+1}}_{j=r_1+1}(d_j-1), 1 \leq i \leq t-1$.}

\medskip

\subsection{The Preinjective to Preinjective Case}

Keep in mind the analysis of $\phi^{II}$ in section 2.3. Let
$$\begin{array}{l} u_1:=\max \{1 \leq j \leq
m^I|c_j>f_1\};\\ v_1:=\max \{1\leq i \leq n^I|c_{u_1+1} \leq f_i\};\\
\cdots\\ u_l:=\max \{1 \leq j \leq m^I|c_j>f_{v_{l-1}+1}\};\\
v_l:=\max \{1\leq i \leq n^I|c_{u_l+1} \leq f_i\};\\
\cdots\\ u_w=m^I. \end{array}$$ Note that $u_1$ is just the number
of the zero blocks in the first block column of $\phi^{II1}$, $v_1$
is just the number of the block columns of $\phi^{II1}$ having the
largest number of nonzero blocks. Dual to Proposition 4 and
Corollary 5 we have:

\medskip

{\bf Proposition 7.} $\rank \phi^{II1}= \min \{\sum^{v_i}_{j=1}f_j+
\sum^{u_w}_{j=u_{i+1}+1}c_j|0 \leq i \leq w-1\}.$

\medskip

{\bf Corollary 8.} $\rank \phi^{II2}=\min \{\sum^{v_i}_{j=1}(f_j-1)+
\sum^{u_w}_{j=u_{i+1}+1}(c_j-1)|0 \leq i \leq w-1\}.$

\medskip

{\bf Remark.} By Proposition 7, Corollary 8, and Proposition 3 we
find that $\rank  \phi^{IIh}, h \in \{1,2\}$, is equal to the rank
of the matrix obtained from $\phi^{IIh}$ by replacing each nonzero
block in $\phi^{IIh}$ with a generic matrix of the same size.

\medskip

{\bf Theorem 9.} {\it $N^I$ is a subrepresentation of $M^I$ if and
only if $v_{w-1}=v_w=n^I, \sum^{v_{w-1}}_{j=v_i+1}f_j \leq
\sum^{u_w}_{j=u_{i+1}+1}c_j$ and $\sum^{v_{w-1}}_{j=v_i+1}(f_j-1)
\leq \sum^{u_w}_{j=u_{i+1}+1}(c_j-1), 0 \leq i \leq w-2$.}

\medskip

\subsection{The Regular to Regular Case}

This case is easier.

\medskip

{\bf Proposition 10.} $\rank  \phi^{RRh}= \sum_{p \in
\mathbb{P}^1(\mathbb{K})} \sum^{\min \{m^p,n^p\}}_{i=1}\min
\{b^p_i,e^p_i\}, h \in \{1,2\}.$

\medskip

{\bf Proof.} Keep in mind the analysis of $\phi^{RR}$ in section
2.4. For every $p \in \mathbb{P}^1(\mathbb{K})$ and every $h \in
\{1,2\}$ we keep the first nonzero element in each row of the
matrices $\phi^{RRhp}_{11}$ and use it to eliminate all other
entries in $\phi^{RRhp}$ which lie in the same row or column by
elementary transformations over the transcendental extension field
$F$ of $\mathbb{K}$ obtained by adding all indeterminates in
$\phi^{RRh}$. Next keep the first nonzero element in each row of the
matrices $\phi^{RRhp}_{22}$, and use them to eliminate all other
elements in $\phi^{RRhp}$ which lie in the same row or column by
elementary transformations over $F$. Proceeding in this way, after
$\min \{m^p,n^p\}$ steps, we obtain $\sum^{\min
\{m^p,n^p\}}_{i=1}\min \{b^p_i,e^p_i\}$ nonzero elements which lie
in different rows and different columns of $\phi^{RRhp}$, while all
other entries in $\phi^{RRhp}$ are reduced to $0$. Thus $\rank
\phi^{RRhp}= \sum^{\min \{m^p,n^p\}}_{i=1}\min \{b^p_i,e^p_i\}.$
Furthermore $\rank \phi^{RRh}= \sum_{p \in \mathbb{P}^1(\mathbb{K})}
\sum^{\min \{m^p,n^p\}}_{i=1}\min \{b^p_i,e^p_i\}, h \in \{1,2\}.$
\hfill{$\Box$}

\medskip

{\bf Theorem 11.} {\it $N^R$ is a subrepresentation of $M^R$ if and
only if $m^p \geq n^p$ and $b^p_i \geq e^p_i, p \in
\mathbb{P}^1(\mathbb{K}), 1 \leq i \leq n^p.$}

\medskip

{\bf Remark.} In the same way, one can show that $N^P$ is a
subrepresentation of $M^I$ if and only if $\sum^{m^I}_{i=1}(c_i-1)
\geq \sum^{n^P}_{i=1}d_i.$ However, to solve the problem completely,
i.e., for arbitrary Kronecker representations $N$ and $M$, more
analysis is needed.

\bigskip

\noindent{\large {\bf Acknowledgement : }} The author would like to
express his gratitude towards the referees for pointing out errors
in earlier versions. He is also grateful to the referees and the
editor for their careful reading and suggestions for improvement.

\medskip

\footnotesize{}

\end{document}